%% file: master.tex
\documentstyle{article}
\newtheorem{theorem}[equation]{Theorem}
\newtheorem{lemma}[equation]{Lemma}
\newtheorem{corollary}[equation]{Corollary}
\newtheorem{proposition}[equation]{Proposition}

\newcommand{\updot}{\textstyle\cdot}

\begin{document}

\input{title}
\input{1}\setcounter{equation}{0}
\input{2}\setcounter{equation}{0}
\input{3}\setcounter{equation}{0}
\input{4}\setcounter{equation}{0}
\input{5}\setcounter{equation}{0}
\input{6}

\input{ref}
\end{document}

%% file: title.tex
\title{Dwork cohomology, de Rham cohomology, \\ and hypergeometric functions}
\author{Alan Adolphson\thanks{Partially supported by NSF grant no.\ 
DMS-9305514.}\\
Department of Mathematics\\
Oklahoma State University\\
Stillwater, Oklahoma 74078\\
adolphs@math.okstate.edu\\
\and
Steven Sperber\\
School of Mathematics\\
University of Minnesota\\
Minneapolis, Minnesota 55455\\
sperber@math.umn.edu} 
\date{}
\maketitle
\begin{center}
{\it Dedicated to Hugh Turrittin on the occasion of his ninetieth birthday}
\end{center}
\begin{abstract}
  In the 1960's, Dwork developed a $p$-adic cohomology theory of de
  Rham type for varieties over finite fields, based on a trace formula
  for the action of a Frobenius operator on certain spaces of
  $p$-analytic functions.  One can consider a purely algebraic
  analogue of Dwork's theory for varieties over a field of
  characteristic zero and ask what is the connection between this
  theory and ordinary de Rham cohomology.  Katz\cite{K1,K2} showed
  that Dwork cohomology coincides with the primitive part of de Rham
  cohomology for smooth projective hypersurfaces, but the exact
  relationship for varieties of higher codimension has been an open
  question.  In this article, we settle the case of smooth, affine,
  complete intersections.
\end{abstract}

%% file: 1.tex
\section{Introduction}

We refer the reader to Katz\cite{K3} for general information on
connections, de Rham cohomology, and the Gauss-Manin connection.  Some
of this material will be reviewed in section 2.  A convenient
reference for the properties of smooth schemes that we use is
Altman-Kleiman\cite[Chap.\ VII, sect.\ 5]{AK}.  Let $S$ be a smooth,
equidimensional ${\bf C}$-scheme and $X$ a smooth, equidimensional
$S$-scheme of relative dimension $N$.  Let $Y\subseteq X$ be a smooth,
closed $S$-subscheme, purely of codimension $r$.  Let
${\mathcal E}$ be a locally free ${\mathcal O}_X$-module of finite
rank with an integrable ${\bf C}$-connection
\[ \nabla:{\mathcal E}\rightarrow\Omega^1_{X/{\bf
    C}}\bigotimes_{{\mathcal O}_X}{\mathcal E}. \]
Let $j:Y\hookrightarrow X$ be the inclusion and let 
\[ \nabla_Y:j^*({\mathcal E})\rightarrow\Omega^1_{Y/{\bf
    C}}\bigotimes_{{\mathcal O}_Y}j^*({\mathcal E}) \] 
be the pullback of $\nabla$ to a connection on $j^*({\mathcal E})$.
We consider the problem of computing the Gauss-Manin connection (an
integrable ${\bf C}$-connection) on the de Rham cohomology sheaves
\[ H^n_{\rm DR}(Y/S,(j^*({\mathcal E}),\nabla_Y)), \]
the hypercohomology of the functor ``direct image under the map
$Y\rightarrow S$'' with respect to the complex
$\Omega^{\updot}_{Y/S}\bigotimes_{{\mathcal O}_Y}j^*({\mathcal E})$. 

In this article, we treat the case where $S,X,Y$ are affine.  Let
$X={\rm Spec}(A)$ and let $f_1,\ldots,f_r\in A$.  Let $y_1,\ldots,y_r$
be indeterminates and consider
\[ {\bf A}_X^r={\rm Spec}(A[y_1,\ldots,y_r]) \]
with projection $\pi:{\bf A}^r_X\rightarrow X$.  Put
\[ F=\sum_{j=1}^r y_jf_j\in A[y_1,\ldots,y_r]. \]
Let $\nabla_{{\bf A}^r_X}$ be the pullback of $\nabla$ to a connection
on $\pi^*({\mathcal E})$.  We let $\nabla_F$ be the ``twist'' of
$\nabla_{{\bf A}^r_X}$ by $\exp F$, i.~e.,
\[ \nabla_F(e)=\nabla_{{\bf A}^r_X}(e)+d_{{\bf A}^r_X/{\bf C}}F\otimes
e, \] where $e$ is a section of $\pi^*({\mathcal E})$ over an open
subset of ${\bf A}^r_X$.  We say that $f_1,\ldots,f_r\in A$ {\it
  define a smooth complete intersection\/} $Y\subseteq X$ if $Y={\rm
  Spec}(A/(f_1,\ldots,f_r))$ is a smooth $S$-scheme of codimension $r$
in $X$.  In particular, $Y$ is regularly immersed in $X$.  

Our main result is the following.
\begin{theorem}
  Suppose $f_1,\ldots,f_r$ define a smooth complete intersection
  $Y\subseteq X$.  Then for every $n\in{\bf Z}$ there is an isomorphism
  of ${\mathcal O}_S$-modules with ${\bf C}$-connections
\[ H^n_{\rm DR}(Y/S,(j^*({\mathcal E}),\nabla_Y))\simeq H^{n+2r}_{\rm
  DR}({\bf A}_X^r/S,(\pi^*({\mathcal E}),\nabla_F)). \] 
\end{theorem}

For example, suppose $S={\rm Spec}(R)$ and $X={\bf A}^N_R$.  Let
${\mathcal E}={\mathcal O}_X$ with the standard connection, so that
$j^*({\mathcal E})={\mathcal O}_Y$ with the standard connection.  Then
the left-hand side is just the de Rham cohomology $H_{\rm DR}^n(Y/S)$ of
the smooth complete intersection $Y\subseteq {\bf A}^N_R$ defined by
$f_1,\ldots,f_r\in R[x_1,\ldots,x_N]$.  We consider the top-dimensional
cohomology group.  Put $C=R[x_1,\ldots,x_N,y_1,\ldots,y_r]$,
the coordinate ring of ${\bf A}^r_X$.  Writing out the connection
$\nabla_F$ in local coordinates, we see that the theorem implies
\begin{equation}
H^{N-r}_{\rm DR}(Y/S)\simeq C\bigg/\biggl(\sum_{i=1}^N
D_{x_i}(C)+\sum_{j=1}^r D_{y_i}(C)\biggr), 
\end{equation}
where
\begin{eqnarray}
D_{x_i} & = & \frac{\partial}{\partial x_i}+\sum_{j=1}^r y_j\frac{\partial
  f_j}{\partial x_i} \\
D_{y_j} & = & \frac{\partial}{\partial y_j}+f_j.
\end{eqnarray}
Furthermore, for $\partial\in{\rm Der}_{\bf C}(R)$, the action of
$\partial$ on $H^{N-r}_{\rm DR}(Y/S)$ via the Gauss-Manin connection is
identified under this isomorphism with the action induced on the
right-hand side of (1.2) by the action of
\begin{equation}
D_{\partial}=\partial+\sum_{j=1}^r y_jf_j^{\partial}
\end{equation}
on $C$, where $\partial$ acts on elements of $C$ (i.~e., polynomials
over $R$) by acting on their coefficients. (In particular,
$f_j^{\partial}$ denotes the polynomial obtained from $f_j$ by
applying $\partial$ to its coefficients.)

In general, by allowing one to replace $Y$ by ${\bf A}^r_X$ and thus
work more globally, the theorem makes the Gauss-Manin connection
easier to compute.  As an application, we show that periods of
differential forms on smooth complete intersections in ${\bf A}^N_R$
satisfy hypergeometric differential equations.  Specifically, we
consider a spanning set of cohomology classes of $H_{\rm
  DR}^{N-r}(Y/S)$ and for each of these cohomology classes we
construct a left ideal in the ring of differential operators on $R$
that annihilates it.  This construction gives a new procedure for
computing Fuchs-Picard equations of complete intersections.

{\it Remark}.  The statement of Theorem 1.1 makes sense only in the
affine case, where one has global regular functions $f_1,\ldots,f_r$
defining the subvariety $Y$ and the connection $\nabla_F$.  However,
even when $X$ is not affine, an analogue of Theorem 4.5 below remains
valid.  In this context, our result is closely related to a theorem of
Hartshorne\cite[Chapter III, Theorem 8.1]{H}.  We plan to return to
this topic in a future article.  The affine case is sufficient for the
application to hypergeometric differential equations treated here. 

The cohomology groups on the right-hand side in Theorem 1.1 are the
``Dwork cohomology'' referred to in the title of this article.  The
first result in the direction of Theorem 1.1 was proved by
Katz\cite{K1,K2}.  Motivated by Dwork's calculations\cite{DW1}, which
showed that the deformation equation in Dwork's $p$-adic cohomology
theory was identical to the corresponding Fuchs-Picard equation, Katz
proved that Dwork cohomology coincides with the primitive part of de
Rham cohomology for smooth projective hypersurfaces.  It was an open
question to determine the relation between these two cohomologies for
smooth complete intersections of codimension greater than one.  This
is what we accomplish here.

This paper is organized as follows.  In section 2, we briefly review
some properties of connections.  In section 3, we use the Leray
spectral sequence of the composition ${\bf A}^r_X\rightarrow
X\rightarrow S$ to reduce Theorem 1.1 to Theorem 4.5.  Theorem 4.5 is
then proved in sections 4 and 5.  In section 6, we show that de Rham
cohomology classes on a complete intersection satisfy hypergeometric
differential equations.

The results of this article have been generalized by Dimca et
al.\cite{DI}.  Using the theory of ${\cal D}$-modules, they obtain
stronger results while avoiding some of the more computational aspects
of our approach.  

%% file: 2.tex
\section{Connections}

We review some basic properties of connections on affine schemes.  Let
$S={\rm Spec}(R)$, where $R$ is a smooth ${\bf C}$-algebra and $X={\rm
  Spec}(A)$, where $A$ is a smooth $R$-algebra.  Let $Y\subseteq X$ be
a smooth, closed $S$-subvariety, say, $Y={\rm Spec}(B)$ with $B=A/I$.
We use algebraic rather than geometric notation, e.~g., we write
$\Omega^1_{A/R}$ rather than $\Omega^1_{X/S}$.  Let $M$ be an
$A$-module with a ${\bf C}$-connection $\nabla$, i.~e., a homomorphism of
abelian groups
\begin{equation}
\nabla:M\rightarrow\Omega^1_{A/{\bf C}}\bigotimes_A M
\end{equation}
satisfying
\[ \nabla(am)=a\nabla(m)+d_{A/{\bf C}}a\otimes m \]
for all $a\in A$, $m\in M$, where $d_{A/{\bf C}}$ denotes the exterior
derivative.  It is also convenient to think of a ${\bf C}$-connection
as an $A$-linear map
\[ \nabla:{\rm Der}_{\bf C}(A)\rightarrow {\rm End}_{\bf C}(M), \]
where ${\rm Der}_{\bf C}(A)$ denotes the Lie algebra of ${\bf
  C}$-linear derivations of $A$.  Let $M_B=M\bigotimes_A B$.  We
define a connection on $M_B$ by pullback.  Specifically, tensoring
(2.1) with $B$ gives a map
\begin{equation}
M_B\rightarrow \Omega^1_{A/{\bf C}}\bigotimes_A M_B\simeq
(\Omega^1_{A/{\bf C}}\bigotimes_A B)\bigotimes_B M_B. 
\end{equation}
The natural map $\Omega^1_{A/{\bf C}}\bigotimes_A B\rightarrow
\Omega^1_{B/{\bf C}}$ induces
\begin{equation}
(\Omega^1_{A/{\bf C}}\bigotimes_A B)\bigotimes_BM_B\rightarrow
\Omega^1_{B/{\bf C}}\bigotimes_BM_B. 
\end{equation}
The composition of (2.2) and (2.3) defines a ${\bf C}$-connection on
the $B$-module $M_B$, which we denote by $\nabla_B$.

When $\nabla$ is an integrable connection, there is an associated de
Rham complex $\Omega^{\updot}_{A/{\bf C}}\bigotimes_A M$ whose
differential we denote by the same symbol used for the connection.
For later use, we recall that
\[ \nabla:\Omega^n_{A/{\bf C}}\bigotimes_A M\rightarrow
\Omega^{n+1}_{A/{\bf C}}\bigotimes_A M \]
is the homomorphism of abelian groups given by
\begin{equation}
\nabla(\omega\otimes m)=d_{A/{\bf C}}\omega\otimes
  m+(-1)^n\omega\wedge\nabla(m), 
\end{equation}
where $\omega\in\Omega^n_{A/{\bf C}}$, $m\in M$, and
$\omega\wedge\nabla(m)$ denotes the image of $\omega\otimes\nabla(m)$
under the canonical map
\[ \Omega^n_{A/{\bf C}}\bigotimes_A\biggl(\Omega^1_{A/{\bf C}}
\bigotimes_A M\biggr)\rightarrow \Omega^{n+1}_{A/{\bf C}}\bigotimes_A
M. \]

A ${\bf C}$-connection on $M$ is, in a natural way, an $R$-connection
via the map $\Omega^1_{A/{\bf C}}\rightarrow\Omega^1_{A/R}$ (or,
equivalently, via the inclusion ${\rm Der}_R(A)\subseteq{\rm Der}_{\bf
  C}(A)$).  By abuse of notation, we also denote the $R$-connection by
$\nabla$.  Furthermore, if $\nabla$ is an integrable connection, then
so is $\nabla_B$.  In this case, $\nabla$ (resp.\ $\nabla_B$) gives
rise to an associated de Rham complex
$(\Omega_{A/R}^{\updot}\bigotimes_A M,\nabla)$ (resp.\ 
$(\Omega_{B/R}^{\updot}\bigotimes_B M_B,\nabla_B)$).  We denote by
$H^n_{\rm DR}(A/R,(M,\nabla))$ (resp.\ $H^n_{\rm
  DR}(B/R,(M_B,\nabla_B))$) the cohomology groups of this complex.  We
wish to compute the $H^n_{\rm DR}(B/R,(M_B,\nabla_B))$ and their
associated Gauss-Manin connections when $M$ is projective of finite
rank.  (We recall that the elements of ${\rm Der}_{\bf C}(R)$ act on
this cohomology group via the Gauss-Manin connection, which is
integrable.)

Let $y_1,\ldots,y_r$ be additional variables and consider ${\bf
  A}^r_A$ with coordinate ring $C=A[y_1,\ldots,y_r]$.  Let
$f_1,\ldots,f_r\in A$ and put
\[ F=\sum_{j=1}^r y_jf_j\in C. \]
Put $M_C=M\bigotimes_A C$ and define a ${\bf C}$-connection $\nabla_C$
on the $C$-module $M_C$ by pullback as before.  We define another
connection $\nabla_F$ on $M_C$ by twisting with $\exp F$.
Specifically, for $\mu\in M_C$, we put
\[ \nabla_F(\mu)=\nabla_C(\mu)+d_{C/{\bf C}}F\otimes\mu. \]
As before, these define $R$-connections as well, and when $\nabla$ is
integrable, so are $\nabla_C$ and $\nabla_F$.  We consider the complex
$(\Omega^{\updot}_{C/R}\bigotimes_C M_C,\nabla_F)$ and its cohomology
groups $H^n_{\rm DR}(C/R,(M_C,\nabla_F))$.

Now suppose that $f_1,\ldots,f_r$ define a smooth complete
intersection $Y\subseteq X$, i.~e., $Y={\rm Spec}(A/I)$,
$I=(f_1,\ldots,f_r)$, is a smooth $R$-scheme of codimension $r$ in
$X$.  Theorem 1.1 can be restated in the following form.
\begin{theorem}
  Suppose $f_1,\ldots,f_r$ define a smooth complete intersection
  $Y\subseteq X$ and $M$ is a projective $A$-module of finite rank.
  Then for every $n\in{\bf Z}$ there is an isomorphism of $R$-modules
  with ${\bf C}$-connections 
\[ H^n_{\rm DR}(B/R,(M_B,\nabla_B))\simeq H^{n+2r}_{\rm
  DR}(C/R,(M_C,\nabla_F)). \] 
\end{theorem}

%% file: 3.tex
\section{Leray spectral sequence}

The purpose of this section is to describe the relation between the
cohomology of the complexes $\Omega^{\updot}_{C/R}\bigotimes_C M_C$
and $\Omega^{\updot}_{C/A}\bigotimes_C M_C$, both with connection
$\nabla_F$.  We denote the cohomology groups of the latter complex by
$H^n_{\rm DR}(C/A,(M_C,\nabla_F))$.  They are $A$-modules with a ${\bf
  C}$-connection (the Gauss-Manin connection), which we denote by
$\delta$.  There is a Leray spectral sequence (\cite[Remark
(3.3)]{K3})
\begin{equation}
E_2^{p,q}=H_{\rm DR}^p(A/R,(H^q_{\rm DR}(C/A,(M_C,\nabla_F)),\delta))
\Rightarrow H^{p+q}_{\rm DR}(C/R,(M_C,\nabla_F)).
\end{equation}

\begin{theorem}
Suppose $f_1,\ldots,f_r$ define a smooth complete intersection
$Y\subseteq X$ and $M$ is a projective $A$-module of finite rank.
Then
\begin{equation}
H^q_{\rm DR}(C/A,(M_C,\nabla_F))=0\qquad\mbox{for $q\neq r$},
\end{equation}
hence the Leray spectral sequence $(3.1)$ collapses and we get
isomorphisms of $R$-modules with ${\bf C}$-connections
\[ H_{\rm DR}^n(A/R,(H^r_{\rm DR}(C/A,(M_C,\nabla_F)),\delta))\simeq 
H^{n+r}_{\rm DR}(C/R,(M_C,\nabla_F))\qquad\mbox{for all $n$.} \]
\end{theorem}

Before beginning the proof, we describe the action of $\nabla_F$ on
$\Omega^{\updot}_{C/A}\bigotimes_C M_C$ in terms of the local
coordinates $y_1,\ldots,y_r$.  An element of
$\Omega^n_{C/A}\bigotimes_C M_C$ can be written uniquely as a sum of
elements of the form
\[ \omega=\mu\,dy_{j_1}\wedge\cdots\wedge dy_{j_n} \]
with $\mu\in M_C$.  Then
\begin{equation}
\nabla_F(\omega)=\biggl(\sum_{j=1}^r D_{y_j}(\mu)\,dy_j\biggr)\wedge
dy_{j_1}\wedge\cdots\wedge dy_{j_n},
\end{equation}
where $D_{y_j}$ is given by (1.4).  Of course, $\mu\in M_C$ can be
uniquely written as a sum of elements of the form $my_1^{a_1}\cdots
y_r^{a_r}$ with $m\in M$, and we have
\begin{equation}
D_{y_j}(my_1^{a_1}\cdots y_r^{a_r})=a_jmy_j^{-1}y_1^{a_1}\cdots
y_r^{a_r}+f_jmy_1^{a_1}\cdots y_r^{a_r}.
\end{equation}

We now define a grading and filtration on $M_C$ that will be used
throughout this article.  An element of $M_C$ can be uniquely written
as a sum of elements of the form $my_1^{a_1}\cdots y_r^{a_r}$ with
$m\in M$.  We define
\[ \deg(m y_1^{a_1}\cdots y_r^{a_r})=a_1+\cdots+a_r. \]
Let $M_C^{(d)}\subseteq M_C$ be the $R$-span of elements $my^a$, $m\in
M$, with $\deg(my^a)=d$.  The corresponding (increasing) filtration
$F.$ on $M_C$ is defined by letting $F_dM_C$ be the $R$-span of
those elements $my^a$ of degree $\leq d$.

{\it Proof of Theorem 3.2}.  We extend the filtration $F.$ to
$\Omega^{\updot}_{C/A}\bigotimes_C M_C$ by defining
\[ F_d(\Omega^{\updot}_{C/A}\bigotimes_C
M_C)=\Omega^{\updot}_{C/A}\bigotimes_C F_dM_C. \] 
To prove the
vanishing result (3.3), it suffices to prove the corresponding result for
the associated graded complex relative to this filtration.  By (3.4) and
(3.5), the associated graded complex is just the Koszul complex on $M_C$
defined by $f_1,\ldots,f_r$.  This Koszul complex obviously decomposes
into a direct sum over $(a_1,\ldots,a_r)\in{\bf N}^r$ of copies of the
Koszul complex on $M$ defined by $f_1,\ldots,f_r$.  Since projective
modules are locally free and cohomology commutes with localization, we
may assume that $M$ is a free $A$-module of finite rank.  We are thus
reduced to proving the following.
\begin{lemma}
Let ${\rm Kos}(A;f_1,\ldots,f_r)$ be the Koszul complex on $A$ defined
by $f_1,\ldots,f_r$.  Then
\[ H^n({\rm Kos}(A;f_1,\ldots,f_r))=0\qquad\mbox{for $n\neq r$.} \]
\end{lemma}

{\it Proof}.  The cohomology groups of this Koszul complex are
$A$-modules and cohomology commutes with localization, so it suffices
to prove the corresponding vanishing of cohomology for each of the
local rings $A_{{\bf p}}$, where ${\bf p}$ is a maximal ideal of $A$.
Since these cohomology groups are supported on $Y$
(\cite[Theorem~16.4]{M}), we may assume ${\bf p}$ corresponds to a
point of $Y$.  But the smooth complete intersection hypothesis implies
that for such ${\bf p}$, the images of $f_1,\ldots,f_r$ in $A_{\bf p}$
form a regular sequence.  The assertion then follows from well-known
properties of regular sequences.

\begin{corollary}
  Suppose $f_1,\ldots,f_r$ define a smooth complete intersection
  $Y\subseteq X$ and $M$ is a projective $A$-module of finite rank.
  Suppose $\mu_j\in F_dM_C$, $j=1,\ldots,r$, are such that
  $\sum_{j=1}^r D_{y_j}(\mu_j)\in F_{d-1}M_C$.  Then there exist
  $\mu'_j\in F_{d-1}M_C$, $j=1,\ldots,r$, such that $\sum_{j=1}^r
  D_{y_j}(\mu'_j)= \sum_{j=1}^r D_{y_j}(\mu_j)$.
\end{corollary}

{\it Proof}.  Let $\mu_j^{(d)}$ denote the homogeneous component of
degree $d$ of $\mu_j$.  The hypothesis implies that $\sum_{j=1}^r
f_j\mu_j^{(d)}=0$, i.~e., $\sum_{j=1}^r
(-1)^{j-1}\mu_j^{(d)}\,dy_1\wedge\cdots\wedge
\widehat{dy_j}\wedge\cdots\wedge dy_r$ is an $(r-1)$-cocycle in
${\rm Kos}(M_C;f_1,\ldots,f_r)$.  But we have just seen that the
$(r-1)$-st cohomology group of this complex vanishes,
hence there exists a skew-symmetric set $\{\eta_{ij}\}_{i,j=1}^r$
(i.~e., $\eta_{ji}=-\eta_{ij}$), with $\eta_{ij}\in M_C^{(d)}$,
such that $\mu_j^{(d)}=\sum_{k=1}^r f_k\eta_{jk}$.  Define
$\mu'_j=\mu_j-\sum_{k=1}^r D_{y_k}(\eta_{jk})$.  Then $\mu'_j\in
F_{d-1}M_C$ and $\sum_{j=1}^r D_{y_j}(\mu'_j)=\sum_{j=1}^r
D_{y_j}(\mu_j)$ by the skew-symmetry of $\{\eta_{ij}\}$.

%% file: 4.tex
\section{The fundamental quasi-isomorphism}

To simplify notation, we put
\[ \bar{M}=H^r_{\rm DR}(C/A,(M_C,\nabla_F)), \]
an $A$-module with the ${\bf C}$-connection $\delta$.  Concretely, by
(3.4) we identify
\begin{equation}
\bar{M}=M_C\bigg/\sum_{j=1}^r D_{y_j}(M_C).
\end{equation}
The connection $\delta$ can be described explicitly as follows.  For
$\mu\in M_C$, let $[\mu]$ denote its image in $\bar{M}$.  Let $m\in M$
and write
\begin{equation}
\nabla(m)=\sum_k\omega_k\otimes m_k,
\end{equation}
with $\omega_k\in\Omega^1_{A/R}$ and $m_k\in M$.  From the definition
of $\nabla_F$ we get
\begin{equation}
\delta([my_1^{a_1}\cdots y_r^{a_r}])=\sum_k
\omega_k\otimes[m_ky_1^{a_1}\cdots y_r^{a_r}] + \sum_{j=1}^r
d_{A/R}f_j\otimes[y_jmy_1^{a_1}\cdots y_r^{a_r}].
\end{equation}
More generally, by (2.4), the map $\delta:\Omega^n_{A/R}\bigotimes_A
\bar{M}\rightarrow \Omega^{n+1}_{A/R}\bigotimes_A \bar{M}$ is given by
\begin{eqnarray}
\delta(\omega\otimes[my_1^{a_1}\cdots y_r^{a_r}]) & = & 
d_{A/R}\omega\otimes[my_1^{a_1}\cdots y_r^{a_r}]+(-1)^n\sum_k
(\omega\wedge\omega_k)\otimes[m_ky_1^{a_1}\cdots y_r^{a_r}] +
\nonumber \\
 & & (-1)^n\sum_{j=1}^r (\omega\wedge d_{A/R}f_j)\otimes[y_jmy_1^{a_1}
\cdots y_r^{a_r}]. 
\end{eqnarray}

Theorem 3.2 gives
\[ H^n_{\rm DR}(A/R,(\bar{M},\delta))\simeq H^{n+r}_{\rm
  DR}(C/R,(M_C,\nabla_F)). \]
In view of this isomorphism, Theorem 2.5 is an immediate consequence
of the following. 
\begin{theorem}
Suppose $f_1,\ldots,f_r$ define a smooth complete intersection
$Y\subseteq X$ and $M$ is a projective $A$-module of finite rank.
Then there is a quasi-isomorphism
\[ \bar{\Phi}:(\Omega^{\updot}_{B/R}\bigotimes_B M_B,\nabla_B)
\rightarrow (\Omega^{\updot}_{A/R}[-r]\bigotimes_A\bar{M},\delta) \]
that induces isomorphisms of $R$-modules with ${\bf C}$-connections
\[ H^n_{\rm DR}(B/R,(M_B,\nabla_B))\simeq H^{n+r}_{\rm
  DR}(A/R,(\bar{M},\delta))\qquad\mbox{for all $n$.} \]
\end{theorem}

The proof of Theorem 4.5 will occupy the remainder of this section and
the next section.  In this section, we define $\bar{\Phi}$, check that
it respects the ${\bf C}$-connections, and show that it is an
isomorphism onto a subcomplex $L^{\updot}[-r]$ of
$\Omega^{\updot}_{A/R}[-r]\bigotimes_A \bar{M}$.  In section 5, we prove
that the inclusion $L^{\updot}\hookrightarrow
\Omega^{\updot}_{A/R}\bigotimes_A \bar{M}$ is a quasi-isomorphism.

We begin by extending the filtration $F.$ on $M_C$ defined in the
previous section to the complex $\Omega^{\updot}_{A/R}\bigotimes_A
\bar{M}$.  Since $\bar{M}$ is a quotient of $M_C$, we get an induced
filtration on $\bar{M}$ which we denote by $F.$ also.  Note that by
Corollary 3.7 there is a natural identification
\begin{equation}
F_d\bar{M}=F_dM_C\bigg/\sum_{j=1}^r D_{y_j}(F_dM_C).
\end{equation}
In what follows, we often make this identification without comment.
We define the filtration $F.$ on $\Omega^{\updot}_{A/R}\bigotimes_A
\bar{M}$ by setting
\[ F_d(\Omega^n_{A/R}\bigotimes_A\bar{M})=\Omega^n_{A/R} \bigotimes_A
F_{d+n-N}\bar{M}, \]
where we denote by $N$ the relative dimension of ${\rm Spec}(A)$ over
${\rm Spec}(R)$.  In particular, $\Omega^1_{A/R}$ is locally free of
rank $N$ and the complex $\Omega^{\updot}_{A/R}\bigotimes_A \bar{M}$
has length $N$.  The indexing is chosen so that $\delta$ respects the
filtration, i.~e., 
\[ \delta(F_d(\Omega^n_{A/R}\bigotimes_A\bar{M}))\subseteq
F_d(\Omega^{n+1}_{A/R}\bigotimes_A\bar{M}). \]

Since $F_d\bar{M}=0$ for $d<0$, we have
\[ F_d(\Omega^n_{A/R}\bigotimes_A\bar{M})=0\qquad\mbox{for $d<N-n$.}
\]
Note that by (4.6),
\[ F_0\bar{M}=F_0M_C\bigg/\sum_{j=1}^r D_{y_j}(F_0M_C). \]
Identifying $F_0M_C$ with $M$ and using $D_{y_j}(M)=f_jM$, we get
\[ F_0\bar{M}=M\bigg/\sum_{j=1}^r f_jM\simeq M_B. \]
Thus for all $n$ we have the identification
\begin{equation}
F_{N-n}(\Omega^n_{A/R}\bigotimes_A\bar{M})=\Omega^n_{A/R}\bigotimes_A
M_B.
\end{equation}

We define a map $\Phi:\Omega^n_{A/R}\bigotimes_A M_B\rightarrow
\Omega^{n+r}_{A/R}\bigotimes_A M_B$ by the formula
\begin{equation}
\Phi(\xi)=(-1)^{nr}d_{A/R}f_1\wedge\cdots\wedge d_{A/R}f_r\wedge\xi,
\end{equation}
where the exterior product on the right-hand side denotes the image of
\[ (d_{A/R}f_1\wedge\cdots\wedge d_{A/R}f_r)\otimes\xi\in
\Omega^r_{A/R}\bigotimes_A\biggl(\Omega^n_{A/R}\bigotimes_A M_B\biggr)
\] 
in $\Omega^{n+r}_{A/R}\bigotimes_A M_B$ under the canonical map.
\begin{lemma}
$\displaystyle \ker\Phi=\sum_{j=1}^r \biggl(\Omega^{n-1}_{A/R}\wedge
d_{A/R}f_j\biggr)\bigotimes_A M_B$
\end{lemma}

{\it Proof}.  It suffices to check equality locally.  Since $M$ is a
projective $A$-module we may assume $M$ is free and thus reduce to the
case $M=A$.  Localizing further if necessary, we get that
$\Omega^1_{A/R}\bigotimes_A B$ is a free $B$-module of rank $N$ and
\[ \Omega^n_{A/R}\bigotimes_A
B\simeq\bigwedge^n\biggl(\Omega^1_{A/R}\bigotimes_A B\biggr) \]
(isomorphism of $B$-modules).  We are thus in the situation of
\cite{S}.  The smooth complete intersection hypothesis implies that
the ideal of $B$ denoted by the symbol ``script-$A$'' in \cite{S} is
the unit ideal.  The desired conclusion then follows from part~(i) of
the theorem of \cite{S}.

Using the identification
\[ \Omega^n_{A/R}\bigotimes_A M_B\bigg/\biggl(\sum_{j=1}^r
\biggl(\Omega^{n-1}_{A/R}\wedge d_{A/R}f_j\biggr)\bigotimes_A
M_B\biggr)\simeq \Omega^n_{B/R}\bigotimes_B M_B \]
and the identification (4.7), we see that, by the lemma, $\Phi$
induces an imbedding 
\[ \bar{\Phi}:\Omega^n_{B/R}\bigotimes_B M_B\hookrightarrow
\Omega^{n+r}_{A/R}\bigotimes_A\bar{M}\qquad\mbox{for all $n$.} \]

We show that $\bar{\Phi}$ is a homomorphism of complexes.  
Every element of $\Omega^n_{B/R}\bigotimes_B M_B$ is a sum of
elements of the form $\bar{\omega}\otimes\bar{m}$ with
$\omega\in\Omega^n_{A/R}$ and $m\in M$, where $\bar{\omega}$ denotes
the image of $\omega$ under the composition
\[ \Omega^n_{A/R}\simeq\Omega^n_{A/R}\bigotimes_A
A\rightarrow\Omega^n_{A/R}\bigotimes_A B\rightarrow\Omega^n_{B/R} \]
and $\bar{m}$ is the image of $m$ under the canonical surjection
$M\rightarrow M_B$.  Thus to show $\bar{\Phi}$ is
a homomorphism of complexes, it suffices to show
\begin{equation}
\delta(\bar{\Phi}(\bar{\omega}\otimes\bar{m}))=\bar{\Phi}(\nabla_B(\bar{\omega}
\otimes\bar{m})). 
\end{equation}
If $\nabla(m)$ is given by (4.2), then the definition of $\nabla_B$
gives
\[ \nabla_B(\bar{m})=\sum_k \bar{\omega}_k\otimes\bar{m}_k, \]
hence by (2.4) we have
\[ \nabla_B(\bar{\omega}\otimes\bar{m}) 
= d_{B/R}\bar{\omega}\otimes\bar{m}+(-1)^n\sum_k
 (\bar{\omega}\wedge\bar{\omega}_k)\otimes\bar{m}_k. \]
This gives
\begin{eqnarray*}
\bar{\Phi}(\nabla_B(\bar{\omega}\otimes\bar{m})) & = & (-1)^{nr+r}
(d_{A/R}f_1\wedge\cdots\wedge d_{A/R}f_r\wedge d_{A/R}\omega)
\otimes\bar{m} + \\
 & & (-1)^{nr+r+n}\sum_k(d_{A/R}f_1\wedge\cdots\wedge
d_{A/R}f_r\wedge\omega\wedge\omega_k) \otimes \bar{m}_k. 
\end{eqnarray*}
From the definition of $\bar{\Phi}$ we have
\[ \bar{\Phi}(\bar{\omega}\otimes\bar{m})=(-1)^{nr}(d_{A/R}f_1\wedge
\cdots\wedge d_{A/R}f_r\wedge\omega)\otimes[m], \] 
hence from (4.4)
\begin{eqnarray*}
\delta(\bar{\Phi}(\bar{\omega}\otimes\bar{m})) & = & (-1)^{nr}d_{A/R}(
d_{A/R}f_1\wedge\cdots\wedge d_{A/R}f_r\wedge\omega)\otimes[m]+ \\
 & & (-1)^{nr+n+r}d_{A/R}f_1\wedge\cdots\wedge d_{A/R}f_r\wedge\omega
\wedge\delta([m]) \\
 & = & (-1)^{nr+r}(d_{A/R}f_1\wedge\cdots\wedge d_{A/R}f_r\wedge
 d_{A/R}\omega)\otimes[m] + \\
 & & (-1)^{nr+n+r}\sum_k (d_{A/R}f_1\wedge\cdots\wedge
 d_{A/R}f_r\wedge\omega\wedge\omega_k)\otimes[m_k].
\end{eqnarray*}
Under the identification (4.7) we have $[m_k]=\bar{m}_k$, which
completes the proof of (4.10).

We now check that the map
\[ H^n_{\rm DR}(B/R,(M_B,\nabla_B))\rightarrow
H^{n+r}_{\rm DR}(A/R,(\bar{M},\delta)) \]
induced by $\bar{\Phi}$ on cohomology respects the Gauss-Manin connection.
We first observe that in the definition of $\bar{\Phi}$, we can replace $R$
by ${\bf C}$.  The map
\[ \Phi_{\bf C}:\Omega^n_{A/{\bf C}}\bigotimes_A M_B\rightarrow
\Omega^{n+r}_{A/{\bf C}}\bigotimes_A M_B \]
defined by
\begin{equation}
\Phi_{\bf C}(\xi)=d_{A/{\bf C}}f_1\wedge\cdots\wedge d_{A/{\bf
    C}}f_r\wedge\xi
\end{equation}
induces a homomorphism of complexes
\[ \bar{\Phi}_{\bf C}:\Omega^{\updot}_{B/{\bf C}}\bigotimes_B
M_B\rightarrow\Omega^{\updot}_{A/{\bf C}}[-r]\bigotimes_A \bar{M}. \]
The Gauss-Manin connection is obtained from the $E_1$-terms of a
spectral sequence associated to a certain filtration on the de Rham
complex over ${\bf C}$.  To check that $\bar{\Phi}$ respects the
Gauss-Manin connection, it suffices to show that $\bar{\Phi}_{\bf C}$
respects this filtration, i.~e., is a homomorphism of filtered
complexes.

These (decreasing) filtrations, which we denote by $G^{\updot}$, are
given by
\begin{eqnarray*}
G^i\Omega^{\updot}_{B/{\bf C}} & = & \mbox{image of }
\biggl(\Omega^i_{R/{\bf C}}\bigotimes_R
B\biggr)\bigotimes_B\Omega^{{\updot}-i}_{B/{\bf C}}\rightarrow
\Omega^{\updot}_{B/{\bf C}} \\
G^i\Omega^{\updot}_{A/{\bf C}} & = & \mbox{image of }
\biggl(\Omega^i_{R/{\bf C}}\bigotimes_R
A\biggr)\bigotimes_{A}\Omega^{{\updot}-i}_{A/{\bf C}}\rightarrow
\Omega^{\updot}_{A/{\bf C}}.
\end{eqnarray*}
The surjection $\Omega^i_{A/{\bf C}}\rightarrow\Omega^i_{B/{\bf C}}$
induces the surjection $\Omega^i_{R/{\bf C}}\bigotimes_R
A\rightarrow\Omega^i_{R/{\bf C}}\bigotimes_R B$.  Thus if
$\bar{\omega}\in G^i\Omega^n_{B/{\bf C}}$, then we can choose a
lifting $\omega\in G^i\Omega^n_{A/{\bf C}}$.  This filtration is
compatible with exterior product, hence
\[ d_{A/{\bf C}}f_1\wedge\cdots\wedge d_{A/{\bf C}}f_r\wedge\omega\in
G^i\Omega^{n+r}_{A/{\bf C}}. \]
It follows from (4.11) that $\bar{\Phi}_{\bf C}(G^i\Omega^{\updot}_{B/{\bf
    C}}\bigotimes_B M_B)\subseteq G^i\Omega^{\updot}_{A/{\bf
    C}}[-r]\bigotimes_A \bar{M}$.

To complete this section, we determine the image of $\bar{\Phi}$.  We
define a subcomplex $L^{\updot}$ of 
$\Omega^{\updot}_{A/R}\bigotimes_A \bar{M}$ by setting
\[ L^n=\{\xi\in F_{N-n}(\Omega^n_{A/R}\bigotimes_A \bar{M})\mid
\delta(\xi)\in F_{N-n-1}(\Omega^{n+1}_{A/R}\bigotimes_A \bar{M})\}. \]
\begin{proposition}
  Suppose $f_1,\ldots,f_r$ define a smooth complete intersection
  $Y\subseteq X$ and $M$ is a projective $A$-module of finite rank.
  Then $\bar{\Phi}$ is an isomorphism of complexes from
  $\Omega^{\updot}_{B/R}\bigotimes_B M_B$ onto $L^{\updot}[-r]$.
\end{proposition}

We begin with a lemma.
\begin{lemma}
Suppose $\omega\in\Omega^{n+r}_{A/R}\bigotimes_A M_B$ satisfies
\[ d_{A/R}f_j\wedge\omega=0\qquad\mbox{for $j=1,\ldots,r$.} \]
Then $\omega\in{\rm im}\;\Phi$.
\end{lemma}

{\it Proof}.  It suffices to check the condition locally, i.~e., to
show that for any maximal ideal ${\bf p}$ of $A$, if
\[ (d_{A/R}f_j)_{\bf p}\wedge\omega_{\bf p}=0\qquad\mbox{for
  $j=1,\ldots,r$,} \] 
then $\omega_{\bf p}\in{\rm im}\;\Phi_{\bf p}$.
Since $\Omega^{\updot}_{A/R}\bigotimes_A M_B$ is supported on $Y$, we
may assume ${\bf p}$ corresponds to a point of~$Y$.  Furthermore,
since $M$ is projective, we may assume that $M$ is a free $A$-module
and thus reduce to the case $M=A$.  But the smooth complete
intersection hypothesis implies that $(d_{A/R}f_1)_{\bf
    p},\ldots,(d_{A/R}f_r)_{\bf p}$ can be extended to a basis of
$(\Omega^1_{A/R})_{\bf p}$ as $A_{\bf p}$-module and $(d_{A/R}f_1)_{\bf
    p}\otimes 1,\ldots,(d_{A/R}f_r)_{\bf p}\otimes 1$ can be
extended to a basis of $(\Omega^1_{A/R})_{\bf p}\bigotimes_{A_{\bf p}}
B_{{\bf p}'}$ as $B_{{\bf p}'}$-module, where ${\bf p}'$ denotes the
image of ${\bf p}$ in $B$.  Since
\[ (\Omega^n_{A/R})_{\bf p}\bigotimes_{A_{\bf p}} B_{\bf p}=
   \bigwedge^n\biggl((\Omega^1_{A/R})_{\bf p}\bigotimes_{A_{\bf p}}
   B_{\bf p}\biggr), \]
the result follows immediately.

{\it Proof of Proposition 4.12}.  We have already proved that
$\bar{\Phi}$ is an injective homomorphism of complexes and (4.10) shows
that its image is contained in $L^{\updot}[-r]$.  So it only remains
to prove $L^{\updot}[-r]\subseteq {\rm im}\;\bar{\Phi}$.  Suppose $\xi\in
F_{N-n}(\Omega^n_{A/R}\bigotimes_A \bar{M})$.  We may write
\[ \xi=\sum_k \omega_k\otimes[m_k], \]
where $\omega_k\in\Omega^n_{A/R}$ and $m_k\in M$.  Write
\[ \nabla(m_k)=\sum_l \omega_{kl}\otimes m_{kl} \]
with $\omega_{kl}\in\Omega^1_{A/R}$ and $m_{kl}\in M$.  By (4.4) we have 
\[ \delta(\xi)=\sum_k d_{A/R}\omega_k\otimes[m_k] + (-1)^n\sum_k \sum_l
(\omega_k\wedge\omega_{kl})\otimes[m_{kl}] + (-1)^n\sum_{j=1}^r 
\sum_k (d_{A/R}f_j\wedge\omega_k)\otimes[y_jm_k]. \]
Using (4.1) and (3.5), we see that $\delta(\xi)\in
F_{N-n-1}(\Omega^{n+1}_{A/R}\bigotimes_A \bar{M})$ (and hence $\xi\in
L^n$) if and only if
\[ d_{A/R}f_j\wedge\biggl(\sum_k\omega_k\otimes[m_k]\biggr)=0\qquad
\mbox{for $j=1,\ldots,r$.} \]
By Lemma 4.13, this implies $\xi\in{\rm im}\;\bar{\Phi}$.

%% file: 5.tex
\section{Completion of the proof}

To complete the proof of Theorem 4.5, it suffices by Proposition 4.12
to show the following.
\begin{theorem}
  Suppose $f_1,\ldots,f_r$ define a smooth complete intersection
  $Y\subseteq X$ and $M$ is a projective $A$-module of finite rank.
  Then the inclusion $L^{\updot}\hookrightarrow
  \Omega^{\updot}_{A/R}\bigotimes_A \bar{M}$ is a quasi-isomorphism.
\end{theorem}

We begin with a lemma.  Let $T$ be a local $R$-algebra that is smooth
over $R$ of relative dimension $N$ and let $f_1,\ldots,f_N$ be
elements of the maximal ideal of $T$ such that
$d_{T/R}f_1,\ldots,d_{T/R}f_N$ form a basis for $\Omega^1_{T/R}$.  Let
$T'$ be a $T$-algebra, let $y_1,\ldots,y_r$ be indeterminates, and consider
\[ \Omega^n_{T/R}\bigotimes_T T'[y_1,\ldots,y_r]. \]
It is a free $T'[y_1,\ldots,y_r]$-module with basis
\begin{equation}
\{d_{T/R}f_{i_1}\wedge\cdots\wedge d_{T/R}f_{i_n}\mid 1\leq
i_1<\cdots<i_n\leq N\}
\end{equation}
and is a free $T'$-module with basis
\begin{equation}
\{y_1^{a_1}\cdots y_r^{a_r}\,d_{T/R}f_{i_1}\wedge\cdots\wedge
d_{T/R}f_{i_n}\mid a_1,\ldots,a_r\in{\bf N},\;1\leq
i_1<\cdots<i_n\leq N\}.
\end{equation}
We grade $T'[y_1,\ldots,y_r]$ by degee, i.~e.,
\[ T'[y_1,\ldots,y_r]^{(d)}=\mbox{$T'$-span of
  $y_1^{a_1}\cdots y_r^{a_r}$ with $a_1+\cdots+a_r=d$} \]
and we grade $\Omega^n_{T/R}\bigotimes_T T'[y_1,\ldots,y_r]$ by
\[ {\rm gr}^{(d)}\biggl(\Omega^n_{T/R}\bigotimes_T
T'[y_1,\ldots,y_r]\biggr)=\Omega^n_{T/R}\bigotimes_T
T'[y_1,\ldots,y_r]^{(d+n-N)}. \]
The map $\Omega^n_{T/R}\bigotimes_T T'[y_1,\ldots,y_r]\rightarrow
  \Omega^{n+1}_{T/R}\bigotimes_T T'[y_1,\ldots,y_r]$ defined by
\begin{equation}
\omega\mapsto\sum_{j=1}^r y_jd_{T/R}f_j\wedge\omega
\end{equation}
then makes $\Omega^{\updot}_{T/R}\bigotimes_T T'[y_1,\ldots,y_r]$ into
a graded complex of $T'$-modules.
\begin{lemma}
With notation and hypotheses as above, we have
\[ H^n({\rm gr}^{(d)}(\Omega^{\updot}_{T/R}\bigotimes_T
T'[y_1,\ldots,y_r]))=0\qquad\mbox{for $d>N-n$.} \]
\end{lemma}

{\it Proof}.  The proof is by induction on $r$.  Suppose $r=1$ and let
\[ \omega\in\Omega^n_{T/R}\bigotimes_T T'[y_1]^{(d+n-N)}. \]
The condition $d>N-n$ implies that $\omega$ is divisible by $y_1$,
i.~e., $\omega$ can be written 
\[ \omega=\sum_{1\leq i_1<\cdots<i_n\leq N}
y_1\omega(i_1,\ldots,i_n)\, d_{T/R}f_{i_1}\wedge\cdots\wedge
d_{T/R}f_{i_n} \]
with $\omega(i_1,\ldots,i_n)\in T'[y_1]^{(d+n-N-1)}$.  The
condition that $\omega$ be a cocycle is that
\[ d_{T/R}f_1\wedge\omega=0. \]
Since (5.2) is a basis, we see that this is the case if and only if
$\omega(i_1,\ldots,i_n)\neq 0$ implies $i_1=1$.  Put
\[ \xi=\sum_{2\leq i_2<\cdots<i_n\leq N} \omega(1,i_2,\ldots,i_n)\,
d_{T/R}f_{i_2}\wedge\cdots\wedge d_{T/R}f_{i_n}. \]
Then $\xi\in{\rm gr}^{(d)}(\Omega^{n-1}_{T/R}\bigotimes_T T'[y_1])$ and
\[ \omega=y_1d_{T/R}f_1\wedge\xi, \]
so $\omega$ is a coboundary.

Now suppose the lemma true for $r-1$ and let
\[ \omega\in{\rm gr}^{(d)}(\Omega^n_{T/R}\bigotimes_T
T'[y_1,\ldots,y_r]). \]
Let $h$ be the highest power of $y_1$ appearing in any term of
$\omega$ (in the decomposition corresponding to the basis (5.3)) and let
$\omega^{(h)}$ be the sum of all terms of $\omega$ of degree $h$ in
$y_1$.  Suppose $h>0$.  Looking at the terms of degree $h+1$ in $y_1$
in the cocycle equation
\begin{equation}
\sum_{j=1}^r y_jd_{T/R}f_j\wedge\omega=0
\end{equation}
gives
\[ d_{T/R}f_1\wedge\omega^{(h)}=0, \]
hence
\[ \omega^{(h)}=\sum_{2\leq i_2<\cdots<i_n\leq N}
y_1\xi(i_2,\ldots,i_n)\, d_{T/R}f_1\wedge
d_{T/R}f_{i_2}\wedge\cdots\wedge d_{T/R}f_{i_n} \]
for some $\xi(i_2,\ldots,i_n)\in T'[y_1,\ldots,y_r]^{(d+n-N-1)}$.  Put
\[ \xi=\sum_{2\leq i_2<\cdots<i_n\leq N}\xi(i_2,\ldots,i_n)\, 
d_{T/R}f_{i_2}\wedge\cdots\wedge d_{T/R}f_{i_n}. \]
Then $\xi\in{\rm gr}^{(d)}(\Omega^{n-1}_{T/R}\bigotimes_T
T'[y_1,\ldots,y_r])$ and the highest power of $y_1$ appearing in
$\omega-(\sum_{j=1}^r y_jd_{T/R}f_j\wedge\xi)$ is $\leq h-1$.

We may thus reduce to the case $h=0$, i.~e., $y_1$ does not appear in
$\omega$.  The cocycle equation (5.6) then implies
\begin{equation}
d_{T/R}f_1\wedge\omega=0
\end{equation}
and
\begin{equation}
\sum_{j=2}^r y_jd_{T/R}f_j\wedge\omega=0.
\end{equation}
From (5.7) we have
\[ \omega=\sum_{2\leq i_2<\cdots<i_n\leq N}\omega(i_2,\ldots,i_n)\,
d_{T/R}f_1\wedge d_{T/R}f_{i_2}\wedge\cdots\wedge d_{T/R}f_{i_n} \]
with $\omega(i_2,\ldots,i_n)\in T'[y_2,\ldots,y_r]^{(d+n-N)}$.
Put
\[ \omega'=\sum_{2\leq i_2<\cdots<i_n\leq N}\omega(i_2,\ldots,i_n)\,
d_{T/R}f_{i_2}\wedge\cdots\wedge d_{T/R}f_{i_n}\in{\rm
  gr}^{(d+1)}(\Omega^{n-1}_{T/R}\bigotimes_T T'[y_2,\ldots,y_r]). \]
From (5.8) we have
\[ \sum_{j=2}^r y_jd_{T/R}f_j\wedge\omega'=0, \]
so by the induction hypothesis there exists
\[ \xi'\in{\rm gr}^{(d+1)}(\Omega^{n-2}_{T/R}\bigotimes_T
T'[y_2,\ldots,y_r]) \]
such that
\begin{equation}
\sum_{j=2}^r y_jd_{T/R}f_j\wedge\xi'=\omega'.
\end{equation}
If we put
\[ \xi=-d_{T/R}f_1\wedge\xi'\in{\rm
  gr}^{(d)}(\Omega^{n-1}_{T/R}\bigotimes_T T'[y_1,\ldots,y_r]), \]
then (5.9) implies
\begin{eqnarray*}
\sum_{j=1}^r y_jd_{T/R}f_j\wedge\xi & = & d_{T/R}f_1\wedge\omega' \\
 & = & \omega,
\end{eqnarray*}
hence $\omega$ is a coboundary.  This completes the proof of Lemma 5.5.

{\it Proof of Theorem 5.1}.  To show that the inclusion
$L^{\updot}\hookrightarrow \Omega^{\updot}_{A/R}\bigotimes_A \bar{M}$
is a quasi-isomorphism, it suffices to show that the corresponding map
of associated graded complexes (relative to the filtration $F.$
defined previously)
\begin{equation}
{\rm gr}^F(L^{\updot})\hookrightarrow{\rm
  gr}^F(\Omega^{\updot}_{A/R}\bigotimes_A \bar{M})
\end{equation}
is a quasi-isomorphism.  But $F.$ induces the ``stupid'' filtration on
$L^{\updot}$:
\[ F_dL^n=\left\{\begin{array}{ll} L^n & \mbox{if $d\geq N-n$,} \\ 
0 & \mbox{if $d<N-n$,} \end{array}\right. \]
hence ${\rm gr}^F_d(L^{\updot})$ is the complex with $L^{N-d}$ in
degree $N-d$ and zeros elsewhere if $0\leq d\leq N$ and is the zero
complex otherwise.  Thus the assertion that (5.10) is a quasi-isomorphism
is equivalent to the assertion that
\[ 0\rightarrow L^{N-d}\rightarrow F_d(\Omega^{N-d}_{A/R}\bigotimes_A
\bar{M})\rightarrow {\rm gr}^F_d(\Omega^{N+1-d}_{A/R}\bigotimes_A
\bar{M})\rightarrow\cdots\rightarrow {\rm
  gr}^F_d(\Omega^N_{A/R}\bigotimes_A \bar{M})\rightarrow 0 \] 
is exact for $0\leq d\leq N$ and
\[ 0\rightarrow {\rm gr}^F_d(\Omega^0_{A/R}\bigotimes_A
\bar{M})\rightarrow\cdots\rightarrow {\rm
  gr}^F_d(\Omega^N_{A/R}\bigotimes_A \bar{M})\rightarrow 0 \] 
is exact for $d>N$.  The definition of $L^{\updot}$ shows that the
sequence 
\[ 0\rightarrow L^{N-d}\rightarrow F_d(\Omega^{N-d}_{A/R}\bigotimes_A
\bar{M})\rightarrow {\rm gr}^F_d(\Omega^{N+1-d}_{A/R}\bigotimes_A
\bar{M}) \]
is exact for $0\leq d\leq N$.  Thus we need to show that the sequence 
\begin{equation}
{\rm gr}^F_d(\Omega^{n-1}_{A/R}\bigotimes_A \bar{M})\rightarrow
{\rm gr}^F_d(\Omega^n_{A/R}\bigotimes_A \bar{M})\rightarrow 
{\rm gr}^F_d(\Omega^{n+1}_{A/R}\bigotimes_A \bar{M}) 
\end{equation}
is exact whenever $d>N-n$, i.~e., that
\[ H^n({\rm gr}^F_d(\Omega^{\updot}_{A/R}\bigotimes_A
\bar{M}))=0\qquad\mbox{whenever $d>N-n$.} \]

There are natural identifications of $A$-modules
\begin{eqnarray*}
{\rm gr}_d^F(\Omega^n_{A/R}\bigotimes_A \bar{M}) & \simeq &
\Omega^n_{A/R}\bigotimes_A \biggl(F_{d+n-N}M_C\bigg/\biggl(F_{d+n-N-1}M_C +
\sum_{j=1}^r D_{y_j}(F_{d+n-N}M_C)\biggr)\biggr) \\
 & \simeq & \Omega^n_{A/R}\bigotimes_A
 \biggl(M_C^{(d+n-N)}\bigg/\sum_{j=1}^r f_jM_C^{(d+n-N)}\biggr) \\
 & \simeq & \Omega^n_{A/R}\bigotimes_A \biggl(M\bigotimes_A
 B[y_1,\ldots,y_r]^{(d+n-N)}\biggr). 
\end{eqnarray*}
We are thus reduced to proving
\begin{equation}
H^n({\rm gr}^{(d)}(\Omega^n_{A/R}\bigotimes_A(M\bigotimes_A
B[y_1,\ldots,y_r])))=0\qquad\mbox{for $d>N-n$,}
\end{equation}
where 
\[ {\rm gr}^{(d)}(\Omega^n_{A/R}\bigotimes_A (M\bigotimes_A
B[y_1,\ldots,y_r])) =
\Omega^n_{A/R}\bigotimes_A (M\bigotimes_A B[y_1,\ldots,y_r]^{(d+n-N)}) \]
and the coboundary map of the complex is given by
\[ \omega\mapsto \sum_{j=1}^r y_jd_{A/R}f_j\wedge\omega. \]

To prove the vanishing result (5.12), we may first localize at a maximal
ideal ${\bf p}$ of $A$.  Since $M$ is a projective $A$-module, we may
assume $M$ is free and thus reduce to the case $M=A$.  Furthermore,
this complex is supported on $Y$, so we may suppose ${\bf p}$
corresponds to a point of $Y$.  We thus want to show that
\[ H^n({\rm gr}^{(d)}((\Omega^n_{A/R})_{\bf p}\bigotimes_{A_{\bf p}}
B_{{\bf p}'}[y_1,\ldots,y_r]))=0\qquad\mbox{for $d>N-n$,} \] 
where ${\bf p}'$ denotes the image of ${\bf p}$ in $B$.  But the
smooth complete intersection hypothesis implies that there exist
$f_{r+1},\ldots,f_N\in A_{\bf p}$ such that $d_{A_{\bf
  p}/R}f_1,\ldots,d_{A_{\bf p}/R}f_N$ form a basis for
$\Omega^1_{A_{\bf p}/R}\simeq(\Omega^1_{A/R})_{\bf p}$.  The result
then follows by applying Lemma 5.5 with $T=A_{\bf p}$ and $T'=B_{{\bf
    p}'}$.

%% file: 6.tex
\section{Hypergeometric equations}

In this section we fix $X={\bf A}^N_R$ and discuss the case of a
smooth complete intersection $Y\subseteq{\bf A}^N_R$ defined by
$f_1,\ldots,f_r\in R[x_1,\ldots,x_N]$.  In the notation of the
previous sections, we have
\begin{eqnarray*}
A & = & R[x_1,\ldots,x_N] \\
B & = & R[x_1,\ldots,x_N]/(f_1,\ldots,f_r) \\
C & = & R[x_1,\ldots,x_N,y_1,\ldots,y_r].
\end{eqnarray*}
We take $M=A$ with the standard connection, hence $M_B=B$ with the
standard connection and we are computing the usual de Rham cohomology
of the variety $Y$.  Theorem 2.5 gives (where $S={\rm Spec}(R)$)
\begin{equation}
H^n_{\rm DR}(Y/S)\simeq H_{\rm DR}^{n+2r}(C/R,(C,\nabla_F)),
\end{equation}
where $\nabla_F$ is the connection
\[ \nabla_F(\mu)=d_{C/R}\mu+d_{C/R}F\otimes\mu \]
with $\mu\in C$, $F=y_1f_1+\cdots+y_rf_r\in C$.  Let $D_R$ be the ring
of ${\bf C}$-linear differential operators with coefficients in $R$,
i.~e., $D_R$ is the subring of ${\rm End}_{\bf C}(R)$ generated by $R$
and ${\rm Der}_{\bf C}(R)$.

We will be mainly interested in $H^{N-r}_{\rm DR}(Y/R)$.  As noted in
the introduction, this is given by (1.2) and the action of a
derivation $\partial\in{\rm Der}_{\bf C}(R)$ on $H^{N-r}_{\rm
  DR}(Y/S)$ is given by (1.5).  We shall see that the differential
equations satisfied by cohomology classes of $(N-r)$-forms on $Y$ can
be determined from the monomials that appear in the $f_j$.  For
$j=1,\ldots,r$, write
\[ f_j=\sum_{i=1}^{\delta_j}\lambda_{j,i}x^{d_{j,i}}, \]
where $\lambda_{j,i}\in R$ and each $d_{j,i}$ is an $N$-tuple of nonnegative 
integers.  Let
\begin{eqnarray*}
F & = & \sum_{j=1}^r y_jf_j \\
 & = & \sum_{j=1}^r\sum_{i=1}^{\delta_j}\lambda_{j,i}x^{d_{j,i}}y^{e_j},
\end{eqnarray*}
where $e_1,\ldots,e_r$ are the standard unit basis vectors in ${\bf
  R}^r$.  We let $[x^uy^v]\in H^{N+r}_{\rm DR}(C/R,(C,\nabla_F))$ be
the cohomology class represented by the differential form
\[ x^uy^v\,dx_1\wedge\cdots\wedge dx_N\wedge dy_1\wedge\cdots\wedge
  dy_r\in \Omega^{N+r}_{C/R}. \]

We define for each $u\in{\bf N}^N$ an $(N-r)$-form
$\omega_u\in\Omega^{N-r}_{B/R}$ as follows.  For any subset
$\sigma=\{i_1,\ldots,i_r\}\subseteq\{1,\ldots,N\}$, put
\[ J_{\sigma}=\det\left[\frac{\partial f_j}{\partial
    x_{i_k}}\right]_{j,k=1,\ldots,r} \] and let $Y_{\sigma}\subseteq
Y$ be the open subset where $J_{\sigma}$ is invertible.  Consider the
differential $(N-r)$-form on $Y_{\sigma}$
\[ \frac{({\rm sgn}\;\sigma)x^u}{J_{\sigma}} dx_1\wedge\cdots\wedge
\widehat{dx_{i_1}}\wedge\cdots\wedge
  \widehat{dx_{i_r}}\wedge\cdots\wedge dx_N, \]
where ${\rm sgn}\;\sigma=\pm 1$ is chosen so that
\[ ({\rm sgn}\;\sigma)\,dx_1\wedge\cdots\wedge
\widehat{dx_{i_1}}\wedge\cdots\wedge
\widehat{dx_{i_r}}\wedge\cdots\wedge dx_N\wedge
dx_{i_1}\wedge\cdots\wedge dx_{i_r}=dx_1\wedge\cdots\wedge dx_N. \]
One checks that these forms agree on overlaps $Y_{\sigma}\cap
Y_{\sigma'}$, hence define a global $(N-r)$-form
$\omega_u\in\Omega^{N-r}_{B/R}$.  Let $[\omega_u]$ denote the
cohomology class of $\omega_u$ in $H^{N-r}_{\rm DR}(Y/S)$.  One checks
that $[\omega_u]$ corresponds to $[x^u]$ under the isomorphism of
Theorem 2.5.  Note that by Theorem 5.1, the cohomology classes $[x^u]$
span $H^{N+r}_{\rm DR}(C/R,(C,\nabla_F))$, hence by Theorem 2.5 the
classes $[\omega_u]$ span $H^{N-r}_{\rm DR}(Y/S)$.

The main purpose of this section is to describe a left ideal
$I_R(u,v)$ of $D_R$ that annihilates $[x^uy^v]$.  By the preceding
remarks, $I_R(u,0)$ then annihilates the cohomology class
$[\omega_u]\in H^{N-r}_{\rm DR}(Y/S)$.  To define these left ideals,
we define them in the ``generic'' case, and take the ``pullback'' to
$Y$.

By the ``generic'' case, we mean the following.  Let $\{\mu_{j,i}\}$
($j=1,\ldots,r$, $i=1,\ldots,\delta_j$) be a collection of
indeterminates and ${\bf C}[\mu]$ the polynomial ring in these
indeterminates.  For $j=1,\ldots,r$, put
\begin{eqnarray*}
f_j^{(\mu)} & = & \sum_{i=1}^{\delta_j}\mu_{j,i}x^{d_{j,i}}, \\
F^{(\mu)} & = & \sum_{j=1}^r y_jf_j^{(\mu)} \\
 & = & \sum_{j=1}^r\sum_{i=1}^{\delta_j}\mu_{j,i}x^{d_{j,i}}y^{e_j}.
\end{eqnarray*}
Let ${\bf A}^N_{{\bf C}[\mu]}$ denote affine $N$-space over ${\bf
  C}[\mu]$ with coordinate ring $A_{\mu}={\bf C}[\mu][x_1,\ldots,x_N]$
and put $C_{\mu}=A_{\mu}[y_1,\ldots,y_r]$, the coordinate ring of
${\bf A}^{N+r}_{{\bf C}[\mu]}$.  Let $\nabla_{F^{(\mu)}}$ be the
connection on $C_{\mu}$ defined by
\[ \nabla_{F^{(\mu)}}(\omega)=d_{C_{\mu}/{\bf C}[\mu]}\omega+
d_{C_{\mu}/{\bf C}[\mu]}F^{(\mu)}\wedge\omega. \] 
We let $[x^uy^v]_{\mu}\in H^{N+r}_{\rm DR}(C_{\mu}/{\bf
  C}[\mu],(C_{\mu},\nabla_{F^{(\mu)}}))$ be the cohomology class
represented by the differential form
\[ x^uy^v\,dx_1\wedge\cdots\wedge dx_N\wedge dy_1\wedge\cdots\wedge dy_r. \]
Let $D_{{\bf C}[\mu]}$ be the ring of differential operators in the
$\partial/\partial\mu_{j,i}$ with polynomial coefficients.  We
describe a collection of differential operators in $D_{{\bf C}[\mu]}$ that
annihilate $[x^uy^v]_{\mu}$.

First we recall the definition of the hypergeometric system associated
to the collection of lattice points $\{(d_{j,i},e_j)\mid
j=1,\ldots,r,\;i=1,\ldots,\delta_j\}\subseteq{\bf R}^{N+r}$ (see
\cite{GKZ}).  Let $E\subseteq{\bf Z}^{\delta_1+\cdots+\delta_r}$ be
the group of relations among these lattice points, i.~e.,
\[ E=\{(b_{j,i})\mid \sum_{j=1}^r\sum_{i=1}^{\delta_j}b_{j,i}(d_{j,i},e_j)=0
\}. \]
For $b=(b_{j,i})\in E$, let $\Box_b$ be the constant coefficient differential
operator
\[ \Box_b=\prod_{b_{j,i}>0}\biggl(\frac{\partial}{\partial\mu_{j,i}}\biggr)
^{b_{j,i}}-\prod_{b_{j,i}<0}\biggl(\frac{\partial}{\partial\mu_{j,i}}\biggr)
^{-b_{j,i}}. \]
Write $d_{j,i}=(d_{j,i}(1),\ldots,d_{j,i}(N))\in{\bf N}^N$.  Define
\begin{eqnarray*}
Z_k & = & \sum_{j=1}^r\sum_{i=1}^{\delta_j} d_{j,i}(k)\mu_{j,i}
\frac{\partial}{\partial\mu_{j,i}} \qquad\mbox{for $k=1,\ldots,N$,} \\
Z_{N+k} & = & \sum_{i=1}^{\delta_k}\mu_{k,i}\frac{\partial}{\partial
\mu_{k,i}} \qquad\mbox{for $k=1,\ldots,r$.}
\end{eqnarray*}
Let $\beta=(\beta_1,\ldots,\beta_{N+r})\in{\bf C}^{N+r}$.  By the {\it
hypergeometric system with parameter $\beta$ associated to the collection of
lattice points} $\{(d_{j,i},e_j)\}_{j,i}$ we mean the system of differential
equations
\begin{eqnarray}
\Box_b(Y) & = & 0 \qquad\mbox{for all $b\in E$,} \\
Z_k(Y) & = & \beta_kY \qquad\mbox{for $k=1,\ldots,N+r$.}
\end{eqnarray}
Recall(\cite{GKZ}, see also \cite{A}) that this system is holonomic.

\begin{theorem}
  The cohomology class $[x^uy^v]_{\mu}\in H^{N+r}_{\rm
    DR}(C_{\mu}/{\bf C}[\mu],(C_{\mu},\nabla_{F^{(\mu)}}))$ satisfies
  the hypergeometric system
\begin{eqnarray*}
\Box_b([x^uy^v]_{\mu}) & = & 0 \qquad\mbox{for $b\in E$,} \\
Z_k([x^uy^v]_{\mu}) & = & -(u_k+1)[x^uy^v]_{\mu} \qquad\mbox{for 
$k=1,\ldots,N$,} \\
Z_{N+k}([x^uy^v]_{\mu}) & = & -(v_k+1)[x^uy^v]_{\mu} \qquad\mbox{for 
$k=1,\ldots,r$,}
\end{eqnarray*}
where we have written $u=(u_1,\ldots,u_N)$, $v=(v_1,\ldots,v_r)$.
\end{theorem}

{\it Proof}.  One computes directly from the definition of
$\nabla_{F^{(\mu)}}$ that there is an isomorphism of $D_{{\bf
    C}[\mu]}$-modules 
\begin{equation}
H^{N+r}_{\rm DR}(C_{\mu}/{\bf C}[\mu],(C_{\mu},\nabla_{F^{(\mu)}}))\simeq 
C_{\mu}\bigg/\biggl(\sum_{i=1}^N D_{x_i}^{(\mu)}(C_{\mu})
+\sum_{j=1}^r D_{y_j}^{(\mu)}(C_{\mu})\biggr),
\end{equation}
where
\begin{eqnarray*}
D_{x_i}^{(\mu)} & = & \frac{\partial}{\partial x_i}+\sum_{j=1}^r y_j
\frac{\partial f_j^{(\mu)}}{\partial x_i}, \\
D_{y_j}^{(\mu)} & = & \frac{\partial}{\partial y_j}+f_j^{(\mu)}.
\end{eqnarray*}
Under this isomorphism, the cohomology class $[x^uy^v]_{\mu}$ in the
left-hand side corresponds to the class of the monomial $x^uy^v$ in
the quotient on the right-hand side.  Furthermore, the action of
$\partial/\partial\mu_{j,i}$ on \linebreak $H^{N+r}_{\rm
  DR}(C_{\mu}/{\bf C}[\mu],(C_{\mu},\nabla_{F^{(\mu)}}))$ is induced by the
action of the differential operator
\begin{eqnarray}
D_{\mu_{j,i}} & = & \frac{\partial}{\partial\mu_{j,i}}+\frac{\partial
F^{(\mu)}}{\partial\mu_{j,i}} \nonumber \\
 & = & \frac{\partial}{\partial\mu_{j,i}}+x^{d_{j,i}}y^{e_j}
\end{eqnarray}
on $C_{\mu}$.  In particular,
\[ \frac{\partial}{\partial\mu_{j,i}}\biggl([x^uy^v]_{\mu}\biggr)=
[x^{u+d_{j,i}}y^{v+e_j}]_{\mu}. \]
More generally, if $\{b_{j,i}\}$ is a collection of nonnegative integers,
\[ \prod_{j,i}\biggl(\frac{\partial}{\partial\mu_{j,i}}\biggr)^{b_{j,i}}
\biggl([x^uy^v]_{\mu}\biggr)=[x^{u+\sum_{j,i}b_{j,i}d_{j,i}}y^{v+\sum_{j,i}
  b_{j,i}e_j}]_{\mu}. \] 
It follows immediately from the definition of
$E$ that $\Box_b([x^uy^v]_{\mu}) =0$ for all $b\in E$.

Using (6.6), we see that for $k=1,\ldots,N$,
\begin{eqnarray*}
Z_k([x^uy^v]_{\mu}) & = & \biggl[\sum_{j=1}^r\sum_{i=1}^{\delta_j}d_{j,i}(k)
\mu_{j,i}D_{\mu_{j,i}}(x^uy^v)\biggr]_{\mu} \\
 & = & \biggl[\sum_{j=1}^r\sum_{i=1}^{\delta_j}d_{j,i}(k)\mu_{j,i}
x^{u+d_{j,i}}y^{v+e_j}\biggr]_{\mu}.
\end{eqnarray*}
Note that 
\begin{eqnarray*}
D_{x_k}^{(\mu)}(x_kx^uy^v) & = & (u_k+1)x^uy^v+\biggl(\sum_{j=1}^r
y_jx_k\frac{\partial f_j^{(\mu)}}{\partial x_k}\biggr)x^uy^v \\
 & = & (u_k+1)x^uy^v+\sum_{j=1}^r\sum_{i=1}^{\delta_j}d_{j,i}(k)\mu_{j,i}
x^{u+d_{j,i}}y^{v+e_j},
\end{eqnarray*}
thus $Z_k([x^uy^v]_{\mu})=-(u_k+1)[x^uy^v]_{\mu}$.  For $k=1,\ldots,r$, we have
\begin{eqnarray*}
Z_{N+k}([x^uy^v]_{\mu}) & = & \biggl[\sum_{i=1}^{\delta_k}\mu_{k,i}
D_{\mu_{k,i}}(x^uy^v)\biggr]_{\mu} \\
 & = & \biggl[\sum_{i=1}^{\delta_k}\mu_{k,i}x^{u+d_{k,i}}y^{v+e_k}\biggr]_{\mu}
\\
 & = & [f_k^{(\mu)}x^uy^{v+e_k}]_{\mu}.
\end{eqnarray*}
But
\[ D_{y_k}^{(\mu)}(x^uy^{v+e_k})=(v_k+1)x^uy^v+f_k^{(\mu)}x^uy^{v+e_k}, \]
thus $Z_{N+k}([x^uy^v]_{\mu})=-(v_k+1)[x^uy^v]_{\mu}$.  This completes
the proof of Theorem 6.4.

We use this result to describe left ideals in $D_R$ that annihilate
the $[x^uy^v]$ and $[\omega_u]$.  Let $I(u,v)\subseteq D_{{\bf
    C}[\mu]}$ be the left ideal generated by the $\Box_b$ for $b\in
E$, the $Z_k+u_k+1$ for $k=1,\ldots,N$, and the $Z_{N+k}+v_k+1$ for
$k=1,\ldots,r$.  By Theorem 6.4, $I(u,v)$ annihilates $[x^uy^v]_{\mu}$.
Let $\phi:{\bf C}[\mu]\rightarrow R$ be the homomorphism of ${\bf
  C}$-algebras defined by setting $\phi(\mu_{j,i})=\lambda_{j,i}$ for
all $j,i$.  Via $\phi$, we regard $R$ as a module over ${\bf C}[\mu]$.
There is a standard construction that associates to any module $M$
over $D_{{\bf C}[\mu]}$ a module $M^{*}$ over $D_R$, namely, take
$M^{*}=R\bigotimes_{{\bf C}[\mu]} M$ with the following action.  We
regard $M^{*}$ as a module over $R$ in the obvious way, while the
action of $\partial\in{\rm Der}_{\bf C}(R)$ on $M^{*}$ is defined by
the formula
\[ \partial(e\otimes m)=\partial(e)\otimes
m+\sum_{j=1}^r\sum_{i=1}^{\delta_j} e\partial(\lambda_{j,i})\otimes\frac
{\partial}{\partial\mu_{j,i}}(m) \]
for $e\in R$.  This construction is functorial, in particular, the
inclusion $\iota:I(u,v)\hookrightarrow D_{{\bf C}[\mu]}$ gives rise to
a homomorphism of $D_R$-modules
\[ \iota^{*}={\rm id.}\otimes\iota:I(u,v)^{*}=R\bigotimes_{{\bf C}[\mu]}I(u,v)
\rightarrow (D_{{\bf C}[\mu]})^{*}=R\bigotimes_{{\bf C}[\mu]}D_{{\bf C}[\mu]}. 
\]
Heuristically, $(D_{{\bf C}[\mu]})^{*}$ consists of formal
differential operators in the $\partial/\partial\mu_{j,i}$ with
coefficients in $R$.  Of course, these differential operators do not
form a ring, in general, but they do form a module over $D_R$ (using
the above definition).  We define a homomorphism $\rho:D_R\rightarrow
(D_{{\bf C}[\mu]})^{*}$ of $D_R$-modules by setting
$\rho(L)=L\cdot(1\otimes 1)$.  Let $I_R(u,v)\subseteq D_R$ be the
inverse image of $\iota^{*}(I(u,v)^{*})$ under $\rho$.  It is a left
ideal of $D_R$.

\begin{theorem}
  The differential operators in the left ideal $I_R(u,v)$ annihilate
  the cohomology class $[x^uy^v]\in H^{N+r}_{\rm
    DR}(C/R,(C,\nabla_F))$.
\end{theorem}

The following corollary is then an immediate consequence of our
discussion of $[\omega_u]$ at the beginning of this section.
\begin{corollary}
  Suppose $f_1,\ldots,f_r$ define a smooth complete intersection
  $Y\subseteq {\bf A}^N_R$.  Then the differential operators in the
  left ideal $I_R(u,0)$ annihilate the cohomology class $[\omega_u]\in
  H^{N-r}_{\rm DR}(Y/S)$.
\end{corollary}

{\it Proof of Theorem 6.7}.  Since $H^{N+r}_{\rm
  DR}(C/R,(C,\nabla_F))$ is given by the right-hand side of (1.2),
  equation (6.5) and the right-exactness of tensor product give
\begin{equation}
H^{N+r}_{\rm DR}(C/R,(C,\nabla_F))\simeq R\bigotimes_{{\bf C}[\mu]}
H^{N+r}_{\rm DR}(C_{\mu}/{\bf C}[\mu],(C_{\mu},\nabla_{F^{(\mu)}}))
\end{equation}
(an isomorphism of $D_R$-modules).  One checks that under this
isomorphism $[x^uy^v]$ corresponds to $1\otimes[x^uy^v]_{\mu}$.

Let $\gamma_{u,v}:D_{{\bf C}[\mu]}\rightarrow H^{N+r}_{\rm
  DR}(C_{\mu}/{\bf C}[\mu],(C_{\mu},\nabla_{F^{(\mu)}}))$ be the
homomorphism of $D_{{\bf C}[\mu]}$-modules defined by setting
$\gamma_{u,v}(L)=L([x^uy^v]_{\mu})$.  By Theorem 6.4, $I(u,v)$ lies in the
kernel of $\gamma_{u,v}$.  Tensoring with $R$, we get a homomorphism
  of $D_R$-modules 
\[ \gamma_{u,v}^{*}:(D_{{\bf C}[\mu]})^{*}\rightarrow
  R\bigotimes_{{\bf C}[\mu]} H^{N+r}_{\rm
  DR}(C_{\mu}/{\bf C}[\mu],(C_{\mu},\nabla_{F^{(\mu)}})) \]
whose kernel contains $\iota^{*}(I(u,v)^{*})$.  Composing with $\rho$,
gives a homomorphism of $D_R$-modules
\[ \gamma_{u,v}^{*}\circ\rho:D_R\rightarrow R\bigotimes_{{\bf C}[\mu]}
  H^{N+r}_{\rm DR}(C_{\mu}/{\bf C}[\mu],(C_{\mu},\nabla_{F^{(\mu)}})) \]
which sends $L$ to $L\cdot(1\otimes[x^uy^v]_{\mu})$ and whose kernel
contains $I_R(u,v)$.  Using the identification~(6.9), we may regard this
as the homomorphism of $D_R$-modules from $D_R$ to $H^{N+r}_{\rm
  DR}(C/R,(C,\nabla_F))$ that sends $L$ to $L([x^uy^v])$.  The fact
that its kernel contains $I_R(u,v)$ is the assertion of the theorem.

{\it Example}.  Consider the Legendre family of (affine) elliptic curves $Y$
defined by $f(x_1,x_2)=0$, where
\[ f(x_1,x_2)=x_2^2-x_1(x_1-1)(x_1-\lambda). \]
We regard this as a smooth affine plane curve over $R={\bf
  C}[\lambda,(\lambda(1-\lambda))^{-1}]$.  The differential
$\omega_0\in\Gamma(Y,\Omega_{Y/S}^1)$ is given on the open set where
$x_2\neq 0$ by $\omega_0=dx_1/2x_2$.  We show that the left ideal
$I_R(0,0)$, which by Corollary~6.8 annihilates $[\omega_0]\in H_{\rm
  DR}^1(Y/S)$, is generated by the Gaussian hypergeometric operator
\begin{equation}
\biggl(\frac{d}{d\lambda}\biggr)^2+\frac{1-2\lambda}{\lambda(1-\lambda)}
\frac{d}{d\lambda}-\frac{1}{4\lambda(1-\lambda)}.
\end{equation}

Multiplying out the formula for $f(x_1,x_2)$, we get
\[ f(x_1,x_2)=x_2^2-x_1^3+(\lambda+1)x_1^2-\lambda x_1. \]
Thus 
\[ f^{(\mu)}(x_1,x_2)=\mu_1x_2^2+\mu_2x_1^3+\mu_3x_1^2+\mu_4x_1 \]
and the map $\phi:{\bf C}[\mu]\rightarrow R$ is defined by
\[ \phi(\mu_1)=1,\quad \phi(\mu_2)=-1,\quad \phi(\mu_3)=\lambda+1,\quad
\phi(\mu_4)=-\lambda. \]
By definition, the left ideal $I(0,0)\subseteq D_{{\bf C}[\mu]}$ is 
generated by the four operators $Z_i+1$ for $i=1,2,3$, where
\begin{eqnarray*}
Z_1 & = & 3\mu_2\frac{\partial}{\partial\mu_2}+2\mu_3\frac{\partial}{\partial
\mu_3}+\mu_4\frac{\partial}{\partial\mu_4} \\
Z_2 & = & 2\mu_1\frac{\partial}{\partial\mu_1} \\
Z_3 & = & \mu_1\frac{\partial}{\partial\mu_1}+\mu_2\frac{\partial}{\partial
\mu_2}+\mu_3\frac{\partial}{\partial\mu_3}+\mu_4\frac{\partial}{\partial\mu_4},
\end{eqnarray*}
and
\begin{equation}
\Box=\biggl(\frac{\partial}{\partial\mu_3}\biggr)^2-\frac{\partial}
{\partial\mu_2}\frac{\partial}{\partial\mu_4}.
\end{equation}

There is a natural map $D_{{\bf C}[\mu]}\rightarrow (D_{{\bf C}[\mu]})^{*}$
which sends $L$ to $1\otimes L$.  To simplify notation, we write $L^{*}$ in
place of $1\otimes L$ and we represent elements of $(D_{{\bf C}[\mu]})^{*}$ as 
differential operators in the $\partial/\partial\mu_i$ with coefficients in 
$R$.  Thus if
\[ L=\sum_k g_k(\mu_1,\mu_2,\mu_3,\mu_4)\biggl(\frac{\partial}{\partial\mu}
\biggr)^k\in D_{{\bf C}[\mu]}, \]
where $k$ is a multi-index, then
\[ L^{*}=\sum_k g_k(1,-1,\lambda+1,-\lambda)\biggl(\frac{\partial}{\partial\mu}
\biggr)^k\in (D_{{\bf C}[\mu]})^{*}. \]
The $D_R$-homomorphism $\rho:D_R\rightarrow
(D_{{\bf C}[\mu]})^{*}$ can then be expressed as follows.  If $P=\sum_{i=1}^n
h_i(\lambda)(\partial/\partial\lambda)^i\in D_R$, then
\begin{equation}
\rho(P)=\sum_{i=1}^n h_i(\lambda)\biggl(\frac{\partial}{\partial\mu_3}-
\frac{\partial}{\partial\mu_4}\biggr)^i\in(D_{{\bf C}[\mu]})^{*}.
\end{equation}
The left ideal $I_R(0,0)\subseteq D_R$ consists of all $P\in D_R$ such
that $\rho(P)\in\iota^{*}(I(0,0)^{*})$.

\begin{lemma}
The $D_R$-module $(D_{{\bf C}[\mu]})^{*}/\iota^{*}(I(0,0)^{*})$ is a
free, rank $2$ $R$-module with basis $1,\partial/\partial \mu_3$.
\end{lemma}

{\it Proof}.  The fact that $Z_i+1\in I(0,0)$ for $i=1,2,3$ allows us
to express $\mu_i\partial/\partial\mu_i$ modulo $I(0,0)$ for $i=1,2,4$
in terms of $\mu_3\partial/\partial\mu_3$.  Explicitly,
\begin{eqnarray}
\mu_1\frac{\partial}{\partial\mu_1} & \equiv & -\frac{1}{2} \quad \mbox{(mod
$I(0,0)$)} \\
\mu_2\frac{\partial}{\partial\mu_2} & \equiv & \mu_4\frac{\partial}{\partial
\mu_4} \quad \mbox{(mod $I(0,0)$)} \\
\mu_4\frac{\partial}{\partial\mu_4} & \equiv &
-\frac{\mu_3}{2}\frac{\partial}{\partial\mu_3}-\frac{1}{4} \quad
\mbox{(mod $I(0,0)$)}. 
\end{eqnarray}
It follows from these equations that, after multiplication by a
monomial in $\mu_1,\mu_2,\mu_4$, every $L\in D_{{\bf C}[\mu]}$ can be
expressed modulo $I(0,0)$ as a polynomial in $\partial/\partial\mu_3$
with coefficients in ${\bf C}[\mu]$.  Since $\lambda$ is invertible in
$R$, this implies that the powers of $\partial/\partial\mu_3$ span
$(D_{{\bf C}[\mu]})^{*}/\iota^{*}(I(0,0)^{*})$ as $R$-module.  From
(6.11) we get
\begin{equation}
(\mu_2\mu_4)\Box=\mu_2\mu_4\biggl(\frac{\partial}{\partial\mu_3}\biggr)^2
-\biggl(\mu_2\frac{\partial}{\partial\mu_2}\biggr)\biggl(\mu_4
\frac{\partial}{\partial\mu_4}\biggr)\in I(0,0).
\end{equation}
Substitution from (6.15) and (6.16) then gives
\begin{equation}
(4\mu_2\mu_4-\mu_3^2)\biggl(\frac{\partial}{\partial\mu_3}\biggr)^2
-2\mu_3\frac{\partial}{\partial\mu_3}-\frac{1}{4}\in I(0,0).
\end{equation}
This implies that
\begin{equation}
\biggl(\frac{\partial}{\partial\mu_3}\biggr)^2+\frac{2(\lambda+1)}
{(\lambda-1)^2}\frac{\partial}{\partial\mu_3}+\frac{1}{4
(\lambda-1)^2}\in \iota^{*}(I(0,0)^{*})
\end{equation}
and also that $1,\partial/\partial\mu_3$ span $(D_{{\bf C}[\mu]})^{*}/
\iota^{*}(I(0,0)^{*})$ as $R$-module.  If there were a nontrivial
$R$-linear relation between $1$ and $\partial/\partial\mu_3$ in this
quotient, it would follow from (6.16) that there is also a nontrivial
$R$-linear relation between $1$ and
$\partial/\partial\mu_3-\partial/\partial\mu_4$.  By (6.12), this
would imply that there exists a nontrivial first-order operator in
$I_R(0,0)$.  Such an operator annihilates $[\omega_0]$, hence there is
a nontrivial $R$-linear relation between $[\omega_0]$ and
$(\partial/\partial\lambda)([\omega_0])$.  But it is well-known (and
easily checked) that
\[ \frac{\partial}{\partial\lambda}([\omega_0])=-\frac{1}{2(\lambda-1)}
[\omega_0]+\frac{1}{2\lambda(\lambda-1)}[\omega_1], \]
giving a nontrivial $R$-linear relation between $[\omega_0]$ and
$[\omega_1]$.  But this contradicts the fact that $[\omega_0]$ and
$[\omega_1]$ form a basis for the de Rham cohomology of the generic
fiber of $Y\rightarrow S$.  This establishes the lemma.

We now consider a second-order operator
\begin{equation}
P=\biggl(\frac{d}{d\lambda}\biggr)^2+A(\lambda)\frac{d}{d\lambda}+B(\lambda)
\in D_R,
\end{equation}
for which we have
\begin{equation}
\rho(P)=\biggl(\frac{\partial}{\partial\mu_3}-\frac{\partial}{\partial\mu_4}
\biggr)^2+A(\lambda)\biggl(\frac{\partial}{\partial\mu_3}-\frac{\partial}
{\partial\mu_4}\biggr)+B(\lambda)\in (D_{{\bf C}[\mu]})^{*}.
\end{equation}
We express $\rho(P)$ modulo $\iota^{*}(I(0,0)^{*})$ in terms of
the basis $1,\partial/\partial\mu_3$.  The procedure is to expand the
expression (6.21) and express $(\partial/\partial\mu_4)^2$,
$(\partial/\partial\mu_3)(\partial/\partial\mu_4)$, and
$\partial/\partial\mu_4$ in terms of $1$, $\partial/\partial\mu_3$,
and $(\partial/\partial\mu_3)^2$.  We can then use (6.19) to express
everything in terms of $1$ and $\partial/\partial\mu_3$.  Using (6.15)
in (6.11) gives
\[ \mu_2\biggl(\frac{\partial}{\partial\mu_3}\biggr)^2-\mu_4\biggl(
\frac{\partial}{\partial\mu_4}\biggr)^2-\frac{\partial}{\partial\mu_4}\in
I(0,0), \]
hence
\begin{equation}
\biggl(\frac{\partial}{\partial\mu_4}\biggr)^2 \equiv
\frac{1}{\lambda}\biggl(\frac{\partial}{\partial\mu_3}\biggr)^2+
\frac{1}{\lambda}\frac{\partial}{\partial\mu_4}\quad(\mbox{mod }
\iota^{*}(I(0,0)^{*})).
\end{equation}
Applying $\partial/\partial\mu_3$ to (6.16) gives 
\[ 2\mu_4\frac{\partial}{\partial\mu_3}\frac{\partial}{\partial\mu_4}+\mu_3
\biggl(\frac{\partial}{\partial\mu_3}\biggr)^2+\frac{3}{2}\frac{\partial}
{\partial\mu_3}\in I(0,0), \]
hence
\begin{equation}
\frac{\partial}{\partial\mu_3}\frac{\partial}{\partial\mu_4} \equiv 
\frac{\lambda+1}{2\lambda}\biggl(\frac{\partial}{\partial\mu_3}\biggr)^2
+\frac{3}{4\lambda}\frac{\partial}{\partial\mu_3}\quad(\mbox{mod }
\iota^{*}(I(0,0)^{*})).
\end{equation}
Equation (6.16) gives
\begin{equation}
\frac{\partial}{\partial\mu_4} \equiv \frac{\lambda+1}{2\lambda}
\frac{\partial}{\partial\mu_3}+\frac{1}{4\lambda}\quad(\mbox{mod }
\iota^{*}(I(0,0)^{*})).
\end{equation}
Expanding (6.21) and substituting (6.22), (6.23), and (6.24) gives
\begin{equation}
\rho(P)\equiv \biggl(\frac{1-2\lambda}{2\lambda^2}+A(\lambda)
\frac{\lambda-1}{2\lambda}\biggr)\frac{\partial}{\partial\mu_3}+ \frac{1}
{4\lambda^2}-A(\lambda)\frac{1}{4\lambda}+B(\lambda) \quad
\mbox{(mod $\iota^{*}(I(0,0)^{*})$)}.
\end{equation}
Solving for $A(\lambda)$ and $B(\lambda)$, we see that $\rho(P)\equiv 0$
(mod $\iota^{*}(I(0,0)^{*})$) if and only if 
\[ A(\lambda)=\frac{1-2\lambda}{\lambda(1-\lambda)} \quad \mbox{and} \quad
B(\lambda)=-\frac{1}{4\lambda(1-\lambda)}. \]
Substituting these expressions in (6.20) gives (6.10).

We have proved that (6.10) is the unique monic second-order operator
in $I_R(0,0)$.  We saw in the proof of Lemma 6.13 that there are no
lower-order operators in $I_R(0,0)$, hence (6.10) generates this left
ideal.